\documentclass[12pt]{article}
\usepackage{amsmath,amssymb,pb-diagram}

\def \1{{\bf 1}}
\def \a{{{\mathfrak a}}}

\def \A{{\mathbb A}}

\def \Aut{{\rm Aut}}
\def \b{{{\mathfrak b}}}

\def \C{{\mathbb C}}

\def \CF{{\cal F}}
\def \CG{{\cal G}}

\def \CO{{\cal O}}

\def \df{\ \begin{array}{c} 
                _{\rm def}\\ ^{\displaystyle =}\end{array}\ }

\def \eps{\varepsilon}

\def \F{{\mathbb F}}

\def \Ga{\Gamma}

\def \GL{{\rm GL}}

\def \Hom{{\rm Hom}}

\def \N{{\mathbb N}}

\def \O{{\cal O}}
\def \p{{{\mathfrak p}}}
\def \ph{\varphi}
\def \prf{\noindent{\bf Proof: }}
\def \P{{\mathbb P}}

\def \Pic{{\rm Pic}}

\def \q{{\mathfrak q}}
\def \qed{\ifmmode\eqno $\square$\else\noproof\vskip 
                12pt plus 3pt minus 9pt \fi}
        \def\noproof{{\unskip\nobreak\hfill\penalty50\hskip2em\hbox{}%
     \nobreak\hfill $\square$\parfillskip=0pt%
     \finalhyphendemerits=0\par}}

\def \Rings{{\rm Rings}}
\def \setminus{\begin{picture}(18,10)\put(4,6)
                {\line(2,-1){10}}\end{picture}}

\def \Spec{{\rm Spec\,}}

\def \tr{{\hspace{1pt}\rm tr\hspace{2pt}}}

\def \Z{\mathbb Z}
\def \={\ =\ }
\def \({\left(}
\def \){\right)}
\def \:{\colon}

\newtheorem{theorem}{Theorem}[section]

\newtheorem{lemma}[theorem]{Lemma}

\newtheorem{proposition}[theorem]{Proposition}

\renewcommand{\setminus}{\hspace{1pt}
    \begin{picture}(10,10)\put(0,5)
    {\line(4,-1){10}}\end{picture}\hspace{1pt}}

\begin{document}

\pagestyle{myheadings} \markright{SCHEMES OVER $\F_1$}

\title{Schemes over $F_1$}
\author{Anton Deitmar}

\date{}
\maketitle
    
\tableofcontents

\newpage
\section*{Introduction}
Jacques Tits wondered in \cite{Tits} if there would exist a ``field of one element'' $\F_1$ such that for a Chevalley group $G$ one has $G(\F_1)=W$, the Weyl group of $G$. 
Recall the Weyl group is defined as $W=N(T)/Z(T)$ where $T$ is a maximal torus, $N(T)$ and $Z(T)$ are the normalizer and the centralizer of $T$ in $G$.
He then showed that one would be able to define a finite geometry with $W$ as automorphism group.

In this paper we will extend the approach of N. Kurokawa, H. Ochiai, and M. Wakayama \cite{KOW} to ``absolute Mathematics'' to define schemes over the field of one element.
The basic idea of the approach of \cite{KOW} is that objects over $\Z$ have a notion of $\Z$-linearity, i.e., additivity, and that the forgetful functor to $\F_1$-objects therefore must forget about additive structures.
A ring $R$ for example is mapped to the multiplicative monoid $(R,\times)$.
On the other hand the theory also requires a base extension functor from objects over $\F_1$ to objects over $\Z$.
Using the analogy of the finite extensions $\F_{1^n}$ as in \cite{Soule}, we are led to define the base extension of a monoid $A$ as
$$
A\otimes_{\F_1}\Z\df \Z[A],
$$
where $\Z[A]$ is the monoidal ring which is defined in the same way as a group ring, using the monoid structure to get a multiplication. 
Based on these two constructions we here lay the foundations of a theory of schemes over $\F_1$.
In \cite{Kato}, K. Kato constructs ``fans'', which are special cases of schemes over $\F_1$. Fans are used to give explicit desingularisations of toric varieties.

This paper was written during a stay at Kyushu University, Japan.
I thank Masato Wakayama and his students for their inspiration and warm hospitality. 
I also thank Zoran Skoda for bringing the paper \cite{Kato} into my attention.

\newpage
\section{Rings over $\F_1$}
In this paper, a ring will always be commutative with unit element.
Recall that a monoid is a set $A$ with an associative composition that has a unit element.
Homomorphisms of monoids are required to preserve units.
In this paper all monoids will be commutative, so $aa'=a'a$ for all $a,a'\in A$.
From now on in the rest of the paper the word `monoid' will always mean `commutative monoid'.
For a monoid $A$ we will write $A^\times$ for the group of invertible elements.

The category of rings over $\F_1$ is by definition the category of monoids.
For a monoid $A$  we will also write $\F_A$ to emphasize that we view $A$ as a ring over $\F_1$.
Let $\F_1\df \{ 1\}$ be the trivial monoid.

For an $\F_1$-ring $\F_A$ we define the \emph{base extension to $\Z$} by
$$
\F_A\otimes\Z\=\F_A\otimes_{\F_1}\Z\= \Z[A].
$$
where $\Z[A]$ is the monoidal ring which is defined like a group ring with the monoid structure of $A$ giving the multiplication.

In the other direction there is the forgetful-functor $F$ which maps a ring $R$ (commutative with one) to its multiplicative monoid $(R,\times)$.

\begin{theorem}\label{adjoint}
The functor of base extension $\cdot\otimes_{\F_1}\Z$ is left adjoint to $F$, i.e., for every ring $R$ and every $\F_A/\F_1$ we have
$$
\Hom_{\rm Rings}(\F_A\otimes_{\F_1}\Z,R)\ \cong\ \Hom_{\F_1}(\F_A,F(R)).
$$
\end{theorem}

\prf
Let $\ph$ be a ring homomorphism from $\F_A\otimes_{\F_1}\Z=\Z[A]$ to $R$. 
Restricting it to $A$ yields a monoid morphism from $A$ to $(R,\times)$.
So we get a map as in the theorem.
Since a ring homomorphism from $\Z[A]$ is uniquely given by the restriction to $A$ this map is injective.
Since on the other hand every monoid morphism from $A$ to $(R,\times)$ extends to a unique ring homomorphism on $\Z[A]$ the claim follows.
\qed

\subsection{Localization}
Let $S$ be a submonoid of $A$.

\begin{lemma}
For a given submonoid $S$ of $A$ there is a monoid $S^{-1}A$ and a monoid homomorphism $\ph$ from $A$ to $S^{-1}A$, determined up to isomorphism with the following property: $\ph(S)\subset \( S^{-1}A\)^\times$ and $\ph$ is universal with this property, i.e., for every monoid $B$, composing with $\ph$ yields an isomorphism
$$
\Hom_{S\to B^\times}(A,B)\cong\Hom(S^{-1}A,B),
$$
where the left hand side describes the set of all monoid homomorphisms $\ph$ from $A$ to $B$ with $\ph(S)\subset B^\times$.
\end{lemma}

\prf
Uniqueness is clear from the universal property. We show existence.
Define $S^{-1} A$ to be the set $A\times S$ modulo the equivalence relation
$$
(m,s)\sim (m',s')\ \Leftrightarrow\ \exists s''\in S : s'' s' m = s'' sm'.
$$
The multiplication in $S^{-1} A$ is given by $(m,s)(m',s')=(mm',ss')$. We also write $\frac ms$ for the element $[(m,s)]$ in $S^{-1}A$.
The map $\ph\colon m\mapsto \frac m1$ has the desired property.
\qed

\subsection{Ideals and spectra}
For two subsets $S,T$ of a monoid $A$ we write $ST$ for the set of all $st$ where $s\in S$ and $t\in T$.
An \emph{ideal} is a subset $\a$ such that $\a A\subset\a$.
For an ideal $\a$ in $A$ the set $\Z[\a]$ is an ideal in $\Z[A]$.
If $\ph\colon \F_A\to \F_B$ is a morphism and $\a$ is an ideal of $B$, then its pre-image $\ph^{-1}(\a)$ is an ideal of $A$.
For a given subset $T$ of $A$ the set $TA$ is the smallest ideal containing $T$.
We call it the ideal \emph{generated by} $T$.

An ideal $\a\ne A$ is called \emph{prime} if $xy\in \a$ implies $x\in\a$ or $y\in\a$.
Equivalently, an ideal $\a$ is prime iff $A\setminus \a$ is a submonoid (compare \cite{Kato}).
We define the \emph{spectrum} of $\F_A$ to be the set $\Spec\F_A$ of all prime ideals in $A$.
Note that this set is never empty, as the set $A\setminus A^\times$ is always a prime ideal.

If $\p$ is a prime ideal, then the set $S_\p=A\setminus\p$ is a submonoid.
We define 
$$
A_\p\df S_\p^{-1}A
$$
to be the \emph{localization at $\p$}.
Note that for the prime ideal $c=A\setminus A^\times$ the natural map $A\to A_c$ is an isomorphism.

We now introduce a topology on $\Spec \F_A$. 
The closed subsets are the empty set and all sets of the form
$$
V(\a)\df \{\p\in\Spec \F_A : \p\supset\a\},
$$
where $\a$ is any ideal.
One checks that $V(\a)\cup V(\b)=V(\a\cap\b)$, and that $\bigcap_{i\in I} V(\a_i)=V\(\bigcup_i\a_i\)$. Thus the axioms of a topology are satisfied.
The point $\eta=\eta_A=\emptyset$ is contained in every non-empty open set.
On the other hand, the point $c=c_A=A\setminus A^\times$ is closed and contained in every non-empty closed set.

\begin{lemma}
For every $f\in A$ the set
$$
V(f)\df \{\p\in\Spec \F_A : f\in\p\}
$$
is closed.
\end{lemma}

\prf 
Let $\a=Af$ be the ideal generated by $f$. Then $V(f)=V(\a)$.
\qed

\section{Schemes over $\F_1$}
\subsection{The structure sheaf}
Let $\F_A$ be a ring over $\F_1$. On the topological space $\Spec \F_A$ we define a sheaf of $\F_1$-rings as follows. 
For an open set $U\subset\Spec \F_A$ we define $\CO(U)$ to be the set of functions \emph{sections} $s\colon U\to \coprod_{\p\in U}A_\p$ such that $s(\p)\in A_\p$ for each $\p\in U$, and such that $s$ is locally a quotient of elements of $A$.
This means that we require for each $\p\in U$ to exist a neighbourhood $V$ of $\p$, contained in $U$, and elements $a,f\in A$ such that for every $\q\in V$ one has $f\notin\q$ and $s(\q)=\frac af$ in $A_\q$.

\begin{proposition}
\begin{enumerate}
\item 
For each $\p\in\Spec \F_A$ the stalk $\CO_\p$ of the structure sheaf is isomorphic to the localization $A_\p$.
\item
$\Ga(\Spec \F_A,\CO)\cong A$.
\end{enumerate}
\end{proposition}

\prf
For (a) define a morphism $\ph$ from $\CO_\p$ to $A_\p$ by sending each element $(s,U_s)$ of $\CO_\p$ to $s(\p)$.
For the injectivity
assume $\ph(s)=\ph(s')$. On some neighbourhood $U$ of $\p$ we have $s(\q)=\frac af$ and $s'(\q)=\frac{a'}{f'}$ for some $a,a',f,f'\in A$. 
This implies that there is $f''\in A$ with $f''\notin\p$ and $f''f'=f''fa'$.
Assume $U$ to be small enough to be contained in the open set
$$
D(f) \= \{\p\in\Spec \F_A : f\notin\p\}.
$$
Then we conclude that $\frac af=\frac {a'}{f'}$ holds in $A_\q$ for every $\q\in U$ and hence $s=s'$.
For the surjectivity let $\frac af\in A_\p$ with $a,f\in A$, $f\notin\p$.
On $U= D(f)$ define a section $s\in \CO(U)$ by $s(\q)=\frac af\in A_\q$.
Then $\ph(s)=s(\p)=\frac af$.
Part (a) is proven.

For part (b) note that the natural map $A\to A_\p$ for $\p\in\Spec \F_A$ induces a map $\psi\colon A\to \Ga(\Spec \F_A,\CO)$.
We want to sow that this map is bijective.
For the injectivity assume $\psi(a)=\psi(a')$.
Then in particular these sections must coincide on the closed point which implies $a=a'$.
For the surjectivity let $s$ be a global section of $\CO$. 
For the closed point $c$ we get an element of $A$ by $a=s(c)\in A_c=A$.
We claim that $s=\psi(a)$.
Since $c$ is contained in every non-empty closed set, its only open neighbourhood is $\Spec \F_A$ itself.
Since $s(c)=a$ we must have $s(\p)=a$ on some neighbourhood of $c$, hence we have it on $\Spec \F_A$, so $s=\psi(a)$.
\qed

\subsection{Monoidal spaces}
A morphism of monoids $\ph\colon A\to B$ is called \emph{local} if $\ph^{-1}(B^\times)=A^\times$.
An \emph{monoidal space} is a topological space $X$ together with a sheaf $\CO_X$ of monoids.
A \emph{morphism of monoidal spaces} $(X,\CO_X)\to (Y,\CO_Y)$ is a pair $(f,f^\#)$, where $f$ is a continuous map $f\colon X\to Y$ and $f^\#$ is a morphism of sheaves $f^\#\colon \CO_Y\to f_*\CO_X$ of monoids on $Y$.
Such a morphism $(f,f^\#)$ is called \emph{local}, if for each $x\in X$ the induced morphism $f_x^\# : \CO_{Y,f(x)}\to \CO_{X,x}$ is local, i.e. satisfies
$$
(f_x^\#)^{-1} \left( \CO_{X,x}^\times\right) \= \CO_{Y,f(x)}^\times.
$$
A \emph{isomorphism of monoidal spaces} is a morphism with a two-sided inverse.
An isomorphism is always local.

\begin{proposition}\label{2.2}
\begin{enumerate}
\item 
For a ring $\F_A$ over $\F_1$ the pair $(\Spec \F_A,\CO_A)$ is a monoidal space.
\item
If $\ph\colon A\to B$ is a morphism of monoids, then $\ph$ induces a morphism of  monoidal spaces
$$
(f,f^\#)\colon \Spec \F_B\to \Spec \F_A,
$$
thus giving a functorial bijection
$$
\Hom(A,B)\ \cong\ \Hom(\Spec \F_B,\Spec \F_A),
$$
where on the right hand side one only admits local morphisms.
\end{enumerate}
\end{proposition} 

\prf
Part (a) is clear.
For (b) suppose we are given a morphism $\ph\colon A\to B$. 
Define a map $f\colon \Spec \F_B\to\Spec \F_A$ that maps $\p$ to the prime ideal $f(\p)=\ph^{-1}(\p)$.
For an ideal $\a$ we have $f^{-1}(V(\a))=V(\ph(\a))$, where $\ph(\a)$ is the ideal generated by the image $\ph(\a)$.
Thus the map $f$ is continuous.
For every $\p\in\Spec \F_B$ we localize $\ph$ to get a morphism $\ph_\p\colon A_{\ph^{-1}(\p)}\to B_\p$.
Since $\ph_\p$ is the localization, it satisfies $\ph_p^{-1}(B_\p^\times)=A_{\ph^{-1}(\p)}^\times$.
For any open set $U\subset \Spec \F_A$ we obtain a morphism
$$
f^\#\colon \CO_A(U)\to \CO_B(f^{-1}(U))
$$
by the definition of $\CO$, composing with the maps $f$ and $\ph$.
This gives a local morphism of local monoidal spaces $(f,f^\#)$.
We have constructed a map
$$
\psi\colon \Hom(A,B)\to \Hom(\Spec \F_B,\Spec \F_A).
$$
We have to show that $\psi$ is bijective.
For injectivity suppose $\psi(\ph)=\psi(\ph')$. For $\p\in\Spec \F_A$ the morphism $f_\p^\#\colon \CO_{A,f(\p)}\to\CO_{B,\p}$ is the natural localization $\ph\colon A_{\ph^{-1}(\p)}\to B_\p$ and this coincides with the localization of $\ph'$.
In particular for $\p=c$ the closed point we get $\ph=\ph'\colon A\to B$.

For surjectivity let $(f,f^\#)$ be a morphism from $\Spec \F_B$ to $\Spec \F_A$.
On global sections the map $f^\#$ gives a monoid morphism 
$$
\ph\colon A=\CO_A(\Spec \F_A)\to f_*\CO_B(\Spec \F_B)=\CO_B(\Spec \F_B)=B.
$$
For every $\p\in\Spec \F_B$ one has an induced morphism on the stalks $\CO_{A,f(\p)}\to \CO_{B,\p}$ or $A_{f(\p)}\to A_\p$ which must be compatible with $\ph$ on global sections, so we have a commutative diagram
$$
\begin{diagram} 
\node{A}\arrow{e,t}{\ph}\arrow{s}
        \node{B}\arrow{s}\\
\node{A_{f(\p)}}\arrow{e,t}{f_\p^\#}
        \node{B_\p}             
\end{diagram}
$$
Since $f(\p)=(f_\p^\#)^{-1}(\p)$ it follows that $f_\p^\#$ is the localization of $\ph$ and hence the claim.
The last part follows since the first bijection preserves isomorphisms by functoriality and an isomorphism on the spectral side preserves closed points and can thus be extended to an isomorphism of the full spectra.
\qed

\subsection{Schemes}
An \emph{affine scheme over $\F_1$} is a monoidal space which is isomorphic to $\Spec \F_A$ for some $A$.

\begin{lemma}\label{affine}
Every open subset of an affine scheme is a union of affine schemes.
\end{lemma}

\prf
Let $U\subset\Spec \F_A$ open. 
Then  there is an ideal $\a$ such that
$$
U\= D(\a)\= \{ \p\in\Spec \F_A : \p \supset\hspace{-11pt}/ \hspace{5pt}\a\}.
$$
So $U$ is the set of all $\p$ such that there exists $f\in \a$ with $f\notin\p$.
For any $f\in A$ let
$$
D(f)\df \{\p\in\Spec \F_A : f\notin\p\}
$$
Then we get
$$
U\=\bigcup_{f\in\a} D(f).
$$
Let $A_f=f^{-1}A = S_f^{-1}A$ where $S_f=\{1,f,f^2,f^3,\dots\}$.
One checks that the open set $D(f)$ is affine, more precisely,
$$
(D(f),\CO_A|_{D(f)})\ \cong\ \Spec \F_{A_f}.
$$
The Lemma follows.
\qed

A monoidal space $X$ is called a \emph{scheme over $\F_1$}, if for every point $x\in X$ there is an open neighbourhood $U\subset X$ such that $(U,\CO_{X}|_U)$ is an affine scheme over $\F_1$.
A \emph{morphism} of schemes over $\F_1$ is a local morphism of monoidal spaces.

A point $\eta$ of a topological space such that $\eta$  is contained in every non-empty open set, is called a \emph{generic point}.

\begin{proposition}
Every connected scheme over $\F_1$ has a unique generic point.
Morphisms on connected schemes map generic points to generic points.
If for an arbitrary scheme $X$ over $\F_1$ we define 
$$
X(\F_1)\df \Hom(\Spec\F_1,X),
$$
then we get
$$
X(\F_1)\ \cong\ \pi_0(X),
$$ 
where the right hand side is the set of connected components of $X$.
\end{proposition}

\prf
The unique generic point of an affine scheme is given by the empty set.
Since an arbitrary scheme is a union of open affine subschemes, it follows that every connected scheme has a unique generic point.
The rest is clear.
\qed

Let $X$ be a scheme over $\F_1$ and let 
$(f,f^\#)$ be a morphism from $X$ to an affine scheme $\Spec \F_A$.
Taking global sections the sheaf morphism $f^\#\colon \CO_A\to f_*\CO_X$ induces a morphism $\ph\: A\to \CO_X(X)$ of monoids.

\begin{proposition}
The map $\psi\: (f,f^\#)\mapsto\ph$ is a bijection
$$
\Hom(X,\Spec \F_A)\ \to\ \Hom(A,\Ga(X,\CO_X)).
$$
\end{proposition}

\prf
Let $\psi(f,f^\#)=\ph$. For each $\p\in X$ one has a local morphism $f_\p^\#\: \CO_{A,f(\p)}\to\CO_{X,\p}$.
Via the map $\CO_X(X)\to\CO_{X,\p}$ the point $\p$ induces an ideal $\tilde\p$ on $\CO_X(X)$ and $f(\p)=\ph^{-1}(\tilde\p)$.
So $f$ is determined by $\ph$.
Further, since $f_\p^\#$ factorizes over $A_{f(\p)}=A_{\tilde\p\circ\ph}\to\CO_X(X)_{\tilde\p}$ it follows that $f_\p^\#(\frac as)=\frac{\ph(a)}{\ph(s)}$ and so $f^\#$ also is determined by $\ph$, so $\psi$ is injective.
For surjectivity reverse the argument.
\qed

The forgetful functor from Rings to Monoids mapping $R$ to $(R,\times)$ extends to a functor
$$
{\rm Schemes}/\Z\ \to\ {\rm Schemes}/\F_1
$$
in the following way: A scheme $X$ over $\Z$ can be written as a union of affine schemes $X=\bigcup_{i\in I}\Spec A_i$ for some rings $A_i$.
We then map it to $\bigcup_{i\in I}\Spec (A_i,\times)$, where we use the gluing maps from $X$.

Likewise the base change $A\mapsto A\otimes_{\F_1}\Z$ extends to a functor
$$
{\rm Schemes}/\F_1\ \to\ {\rm Schemes}/\Z
$$
by writing a scheme $X$ over $\F_1$ as a union of affine ones, $X=\bigcup_{i\in I}\Spec A_i$ and then mapping it to $\bigcup_{i\in I}\Spec (A_i\otimes_{\F_1}\Z)$, which is glued via the gluing maps from $X$.
The fact that these constructions do not depend on the choices of affine coverings follows from Proposition \ref{2.2}, Lemma \ref{affine} and its counterpart for schemes over $\Z$.

As an example for a scheme which is not affine let us construct the projectiove line $\P^1$ over $\F_1$.
Let $C_\infty=\{ \dots,\tau^{-1},1,\tau,\dots\}$ be the infinite cyclic group with generator $\tau$.
Let $\C_{\infty,+}=\{ 1,\tau,\tau^2,\dots\}$ and $C_{\infty,-}=\{ 1,\tau^{-1},\tau^{-2},\dots\}$.
The inclusions give maps from $U=\Spec C_\infty$ to $X=\Spec C_{\infty,+}$ and $Y=\Spec C_{\infty,-}$ identifying $U$ with open subsets of the latter.
We define a new space $\P^1$ by gluing $X$ and $Y$ along this common open subset.
The space $X$ has two points, $c_X,\eta_X$, one closed and one open and likewise for $Y$.
The space $\P^1$ has three points, $c_X, c_Y, \eta$, two closed and one open.
The structure sheaves of $X,Y,U$ give a structure sheaf of $\P^1$ making it a scheme over $\F_1$.

\section{Fibre products}
Let $S$ be a scheme over $\F_1$. A pair $(X,f_X)$ consisting of an $\F_1$-scheme $X$ and a morphism $f_X\colon X\to S$ is called a \emph{scheme over $S$}

\begin{proposition}
Let $X,Y$ be schemes over $S$. There exists a scheme $X\times_S Y$ over $S$, unique up to $S$-isomorphism and morphisms from $X\times_S Y$ to $X$ and $Y$ such that the diagram
$$
\begin{diagram} 
\node{X\times_SY} \arrow{e} \arrow{s} \arrow{se}
        \node{X} \arrow{s,r}{f_X}\\
\node{Y} \arrow{e,b}{f_Y}
        \node{S}
\end{diagram}
$$
is commutative and the composition with these morphisms induces a bijection for every scheme $Z$ over $S$,
$$
\Hom_S(Z,X)\times\Hom_S(Z,Y)\ \to\ \Hom_S(Z,X\times_SY).
$$
This fibre product is compatible with $\Z$ extension and the usual fibre product of schemes, i.e., one has
$$
\( X\times_S Y\)\otimes\Z\ \cong\ (X\otimes\Z) \times_{S\otimes Z} (Y\otimes\Z).
$$
\end{proposition}

\prf
This Proposition follows via the gluing procedure once it has been established in the affine case.
So let's assume that $X,Y,S$ are all affine, say $X=\Spec A$, $Y=\Spec B$ and $S=\Spec L$.
Then $f_X$ and $f_Y$ give monoid morphisms $\ph_A\colon L\to A$ and $\ph_B\colon L\to B$ and we define
$$
X\times_S Y\df \Spec\( A\otimes_L B\),
$$
where $A\otimes_L B$ is the monoid $A\times B/\sim$ and $\sim$ is the equivalence relation generated by $(\ph_A(l)m,n)\sim(m,\ph_B(l)n)$ for all $m\in A$, $n\in B$, and $l\in L$.
In other  words, $A\otimes_L B$ is the push-out of $\ph_A$ and $\ph_B$ and so $X\times_S Y$ has the desired properties.
\qed

\section{$\CO_X$-modules}
For an $\F_1$-ring $\F_A$, one defines an \emph{$\F_A$-module} to be a set $M$ together with an action $A\times M\to M$, $(a,m)\mapsto am$ such that $1s=s$ and $ (aa')s=a(a's)$ for every $s\in S$ and all $a,a'\in A$.
One defines $M\otimes\Z=\Z[M]$ with the obvious $\F_A\otimes\Z$-module structure.
The direct sum $M\oplus N$ of $A$-modules is the disjoint union and the tensor product $M\otimes_A N$ is the quotient $M\times N /\sim$, where $\sim$ is the equivalence relation generated by $(am,n)\sim (m,an)$ for all $a\in A, m\in M, n\in N$.
Then $M\otimes_A N$ is an $A$-module via $a[m,n]=[am,n]$.
There are natural isomorphisms of $\F_A\otimes\Z$-modules
\begin{eqnarray*}
(M\oplus N)\otimes\Z & \cong& (M\otimes\Z)\oplus (N\otimes\Z)\\
(M\otimes_A N)\otimes\Z &\cong& (M\otimes\Z)\otimes_{A\otimes \Z}(N\otimes\Z). 
\end{eqnarray*}

Let $(X,\CO_X)$ be a monoidal space.
We define an \emph{$\CO_X$-module} to be a sheaf $\CF$ of sets on $X$ together with the structure of an $\CO_X(U)$-module on $\CF(U)$ for each open set $U\subset X$ such that for open sets $V\subset U\subset X$ the restriction $\CF(U)\to\CF(V)$ is compatible with the module structure via the map $\CO_X(U)\to\CO_X(V)$.
If $\CF$ is an $\CO_X$-module and $U\subset X$ is open, then $\CF|_U$ is an $\CO_U=\CO_X|_U$-module.
A \emph{morphism} of $\CO_X$-modules $\ph\colon\CF\to\CG$ is a morphism of sheaves such that for every open set $U\subset X$ the map $\ph(U)\colon\CF(U)\to\CG(U)$ is a morphism of $\CO_X(U)$-modules.
The category of $\CO_X$-modules has kernels, cokernels, images, and internal Hom's, as the pre-sheaf $U\mapsto \Hom_{\CO_X(U)}(\CF(U),\CG(U))$ is a sheaf called the \emph{sheaf Hom} for any two given $\O_X$-modules $\CF$ and $\CG$. 
This then is naturally an $\CO_X$-module.
The \emph{tensor product} $\CF\otimes_{\CO_X}\CG$ of two $\CO_X$-modules $\CF,\CG$ is the sheaf associated to the pre-sheaf $U\mapsto \CF(U)\otimes_{\CO_X(U)}\CG(U)$.
Note that $\CF\otimes\CO_X\cong\CF$. 
An $\CO_X$-module is \emph{free} if it is isomorphic to a direct sum of copies of $\CO_X$.
It is \emph{locally free} if every $x\in X$ has an open neighbourhood $U$ such that $\CF|_U$ is free.
In this case the \emph{rank} of $\CF$ at a point $x$ is the number of copies of $\CO_X$ needed over any open neighbourhood of $x$ over which $\CF$ is free.
If $X$ is connected, then the rank is the same everywhere.

For a given $\CO_X$-module $\CF$ let $\CF^*\df \Hom(\CF,\CO_X)$ be its \emph{dual} module.
There is a canonical morphism, called the \emph{trace},
$$
\tr\colon \CF\otimes\CF^*\ \to\ \CO_X,
$$
given over an open set $U$ by mapping $f\otimes\alpha\in\CF(U)\otimes_{\CO_X(U)}\CF^*(U)$ to $\alpha(f)\in\CO_X(U)$.

\begin{lemma}
If $\CF$ is locally free of rank one, then the trace is an isomorphism.
So the set of isomorphism classes of locally free $\CO_X$-modules of rank one forms a group, called the \emph{Picard group} $\Pic(X)$ of $X$.
\end{lemma}

\prf
Clear.
\qed

\begin{proposition}
Let $X$ be an affine scheme over $\F_1$, then every locally free sheaf is free, so in particular, $\Pic(X)$ is trivial.
\end{proposition}

\prf
There is a unique closed point $c$ which is contained in every non-empty closed set, so the only open neighbourhood of $c$ is the full space $X$, which implies that every locally free sheaf is free.
\qed

\begin{proposition}
The Picard group of the projective line $\P^1$ is isomorphic to $\Z$.
\end{proposition}

\prf
The space $\P^1$ has three elements, $c_X$, $c_Y$, and $\eta$.
The non-trivial open sets are $U=\{\eta\}$, $X=\{\eta, c_X\}$, and $Y=\{\eta, c_y\}$.
We have $\CO(X)\cong C_{\infty,+}$, $\CO(Y)\cong C_{\infty,-}$ and $\CO(U)\cong C_\infty$ with the inclusions as restriction maps.
Since $X$ and $Y$ are affine, a given invertible sheaf $\CF$ is trivial on $X$ and on $Y$.
We fix isomorphisms $\alpha\colon\CF|_X\to \CO_X$ and $\beta\colon\CF|_Y\to\CO_Y$.
The restriction of $\alpha$ and $\beta$ gives two isomorphisms $\CF(U)\to \CO(U)=C_\infty$.
These two differ by a $C_\infty$-module automorphism of $C_\infty$.
The group of these automorphisms is isomorphic to $\Z$.
It is easy to see that this establishes the claimed isomorphism.
\qed

Let $f\colon X\to Y$ be a morphism of $\F_1$-ringed spaces.
If $\CF$ is an $\CO_X$-module, then $f_*\CF$ is an $f_*\CO_X$-module.
The morphism
$f^\#\colon\CO_Y\to f_*\CO_X$ thus makes $f_*\CF$ into an $\CO_Y$-module, called the \emph{direct image} of $\CF$.

For an $\CO_Y$-module $\CG$ the sheaf $f^{-1}\CG$ is an $f^{-1}\CO_Y$-module.
Recall that the functor $f^{-1}$ is adjoint to $f_*$, this implies that
$$
\Hom_X(f^{-1}\CO_Y,\CO_X)\ \cong\ \Hom_Y(\CO_Y,f_*\CO_X).
$$
So the map $f^\#\colon \CO_Y\to f_*\CO_X$ gives a map $f^{-1}\CO_Y\to\CO_X$.
We define $f^*\CG$ to be the tensor product
$$
f^{-1}\CG\otimes_{f^{-1}\CO_Y}\CO_X.
$$
So $f^*\CG$ is an $\CO_X$-module, called the \emph{inverse image} of $\CG$.
The functors, $f_*$ and $f^*$ are adjoint in the sense that
$$
\Hom_{\CO_X}(f^*\CG,\CF)\ \cong\ \Hom_{\CO_Y}(\CG,f_*\CF).
$$

\subsection{Localization}
Let $M$ be a module of the monoid $A$.
Let $S\subset A$ be a submonoid.
We define the localization $S^{-1}M$ to be the following module of $S^{-1}A$.
As a set, $S^{-1}M$ is the set of all pairs $(m,s)=\frac ms$ with $m\in M$ and $s\in S$ modulo the equivalence relation
$$
\frac ms\sim\frac {m'}{s'}\ \Leftrightarrow \ \exists s''\in S : s'' s' m= s'' s m'.
$$
The $S^{-1}A$-module structure is given by
$$
\frac as\cdot \frac m{s'}\= \frac{am}{ss'}.
$$
A given $A$-module homomorphism $\ph\colon M\to N$ induces an $S^{-1}A$-module homomorphism $S^{-1}\ph\colon S^{-1}M\to S^{-1} N$ by $S^{-1}\ph(\frac ms)=\frac{\ph(m)}s$.
Note that $S^{-1}\ph$ is injective/surjective if $\ph$ is.

Given an $A$-module $M$ we define an $\CO_{\F_A}$-sheaf $\tilde M$ on $X=\Spec\F_A$ as follows.
For each prime ideal $\p$ of $A$ let $M_\p$ be the localization $S_\p^{-1} M$ at $\p$.
For any open set $U\subset \Spec\F_A$ we define the set $\tilde M(U)$ to be the set of functions $s\colon U\to\amalg_{\p\in U} M_\p$ such that $s(\p)\in M_\p$  for each $\p$ and such that $s$ is locally a fraction, i.e., we require that for each $\p\in U$ there is a neighbourhood $V\subset U$ of $\p$ and $m\in M$ as well as $f\in A$ with $f\notin\q$ for every $q\in V$ and $s(\q)=\frac mf$ in $m_\q$.
Then $\tilde M$ is a sheaf with the obvious restriction maps.
It is an $\CO_X$-module.
For each $\p\in X$ the stalk $(\tilde M)_\p$ coincides with the localization $M_\p$ at $\p$.
For every $f\in M$ the $A_f$-module $\tilde A(D(f))$ is isomorphic to the localized module $M_f$.
In particular, we have $\Ga(X,\tilde M)\cong M$.
One  also has $(M\otimes_A N)^\sim\cong \tilde M\times_{\CO_X}\tilde N$
For a morphism $f\colon \Spec B\to\Spec A$ one has $f_*\tilde M\cong (_AM)^\sim$, where $_AM$ is $M$ considered as $A$-module via the map $A\to B$ induced by $f$.
Finally, we have $f^*(\tilde N)\cong (B\otimes_A N)^\sim$.

\subsection{Coherent modules}
Let $(X,\CO_X)$ be a scheme over $\F_1$.
An $\CO_X$-module $\CF$ is called \emph{coherent} if every $x\in X$ has an affine neighbourhood $U\cong\Spec\F_A$ such that $\CF\_U$ is isomorphic to $\tilde M$ for some $A$-module $M$.

\begin{proposition}
Let $X$ be a scheme over $\F_1$.
An $\CO_X$-module $\CF$ is coherent if and only if for every open affine subset $U=\Spec A$ of $X$ there is an $A$-module $M$ such that $\CF|_U\cong \tilde M$.
\end{proposition}

\prf
Let $\CF$ be coherent and let $U=\Spec A$ be an affine open subset.
Let $X=\bigcup_i\Spec A_i$ be a covering by affines such that $\CF|_{\Spec A_i}\cong \tilde M)i$ for some modules $M_i$.
Let $f\in A_i$. Then $\CF|_{D(f)}\cong(M_{i,f})^\sim$.
So $X$ has a basis of the topology consisting of affines on which $\CF$ comes from modules.
It follows that $\CF|_U$ also is coherent, so we can reduce to the case when $X$ is affine.
Then $\CF$ must come from a module in a neighbourhood of the closed point, so $\CF$ comes from a module.
\qed

\section{Chevalley Groups}
\subsection{$\GL_n$}
We first repeat the definition of $\GL_n(\F_1)$ as in \cite{KOW}.
On rings the functor $\GL_n$, which is a representable group functor, maps $R$ to $\GL_n(R)=\Aut_R(R^n)$.
To define $\GL_n(\F_1)$ we therefore have to define $\F_1^n$, or more generally, modules over $\F_1$.
Since $\Z$-modules are just additive abelian groups by forgetting the additive structure one simply gets sets.
So $\F_1$ modules are sets.
For an $\F_1$-module $S$ and a ring $R$ the base extension is the free $R$-module generated by $S$: $S\otimes_{\F_1} R = R[S] = R^{(S)}$.
In particular, $\F_1^n\otimes\Z$ should be isomorphic to $\Z^n$, so $\F_1^n$ is a set of $n$ elements, say $\F_1^n=\{ 1,2,\dots,n\}$.
Hence
\begin{eqnarray*}
\GL_n(\F_1) &=& \Aut_{\F_1}(\F_1^n)\\
&=& \Aut(1,\dots,n)\\
&=& S_n
\end{eqnarray*}
the symmetric group in $n$ letters, which happens to be the Weyl group of $\GL_n$.
We now extend this to rings over $\F_1$.
One would expect that $\F_A^n=\F_1^n\otimes \F_A=\{1,\dots,n\} \times A$.
So we define an $\F_A$-module to be a set with an action of the monoid $A$ and in particular $\F_A^n$ is a disjoint union of $n$-copies of $A$.
We define
$$
\GL_n(\F_A)\df \Aut_{\F_A}(\F_A^n).
$$
This is compatible with $\Z$-extension,
\begin{eqnarray*}
\GL_n(\F_A\otimes\Z) &=& \Aut_{\F_A\otimes\Z}(\F_A^n\otimes\Z)\\
&=& \Aut_{\Z[A]}(\Z[A]^n)\\
&=& \GL_n(\Z[A])
\end{eqnarray*}
as required.

Note that $\GL_n(\F_A)$ can be identified with the group of all $n\times n$ matrices $A$ over $\Z[A]$ with exactly one non-zero entry in each row and each column and this entry being in the group of invertible elements $A^\times$.

\begin{proposition}
The group functor $\GL_1$ on $\Rings/\F_1$ is represented by the infinite cyclic group $C_\infty$.
This is compatible with $\Z$-extension as
$$
\F_{C_\infty}\otimes\Z\=\Z[C_\infty]\ \cong\ \Z[T,T^{-1}]
$$
represents $\GL_1$ on rings.
\end{proposition}

\prf
Choose a generator $\tau$ of $C_\infty$. For any ring $\F_A$ over $\F_1$ we have an isomorphism
$$
\Hom(\F_{C_\infty},\F_A)\ \to\ A^\times =\GL_1(\F_A)
$$
mapping $\alpha$ to $\alpha(\tau)$.
\qed

The functor $\GL_n$ on rings over $\F_1$ cannot be represented by a ring $\F_A$ over $\F_1$ since $\Hom(\F_A,\F_1)$ has only one element.

\begin{proposition}
The functor $\GL_n$ on rings over $\F_1$ is represented by a scheme over $\F_1$.
\end{proposition}

\prf
We will give the proof for $\GL_2$. 
The general case will be clear from that.
For $\F_A/\F_1$ the group $\GL_2(\F_A)$ can be identified with the group of matrices
$$
\left(\begin{array}{cc}A^\times &  \\ & A^\times\end{array}\right)
\ \cup\ \left(\begin{array}{cc} & A^\times \\A^\times & \end{array}\right).
$$
We define a scheme $X$ over $\F_1$ as the disjoint union $X_0\cup X_1$ and $X_0\cong X_1\cong \Spec((C_\infty\times C_\infty)$, where $C_\infty$ is the infinite cyclic group.
The group structure on $\Hom(\Spec\F_A,X)$ for every $\F_A/\F_1$ comes via a multiplication map $m\colon X\times_{\F_1} X\to X$ defined in the following way.
The scheme $X\times_{\F_1} X$ has connected components $X_i\times_{\F_1}X_j$ for $i,j\in\{ 0,1\}$. 
The multiplication map splits into components $m_{i,j}\colon X_i\times_{\F_1}X_j\to X_{i+j\ (2)}$.
Each $m_{i,j}$ in turn is given by a monoid morphism $\mu_{i,j}\colon C_\infty^2\to C_\infty^2\times C_\infty^2$, called the comultiplication.
Here $\mu_{i,j}$ maps $a$ to $(\eps^i(a),\eps^j(a))$, where $\eps^0(a)=a$ and $\eps^1(x,y)=(y,x)$. 
This finishes the construction and the proof.
\qed

\subsection{${\rm O}_n$ and ${\rm Sp}_{2n}$}
The orthogonal group ${\text O}_n$ is the subgroup of $\GL_n$ consisting of all matrices $A$ with $AqA^t=q$, where $q={\rm diag}(J,\dots,J,1)$, and $J=\left(\begin{array}{cc} & 1 \\1 & \end{array}\right)$. 
The last 1 only occurring if $n$ is odd.
A computation shows that
$$
\text{O}_n(\F_1)\ \cong\ \text{Weyl\ group}
$$
holds here as well.
Finally $\text{Sp}_{2l}$ is the group of all $A$ with $ASA^t=S$, where $S$ is the $2l\times 2l$ matrix with anti-diagonal $(1,\dots,1,-1,\dots,-1)$ and zero elsewhere.
Likewise, $\text{Sp}_2l(\F_1)$ is the Weyl group.
Both $\text O_n$ and $\text{Sp}_{2l}$ are representable by $\F_1$-schemes.

\section{Zeta functions}
Let $X$ be a scheme over $\F_1$.
For $\F_A/\F_1$ we write as usual,
$$
X(\F_A)\=\Hom(\Spec \F_A,X)
$$ 
for the set of $\F_A$-valued points of $X$.
After Weil we set for a prime $p$,
$$
Z_X(p,T)\df \exp\(\sum_{n=1}^\infty\frac{T^n}n \#(X(\F_{p^n}))\),
$$
where, of course $\F_{p^n}$ means the field of $p^n$ elements and $X(\F_{p^n})$ stands for $X((\F_{p^n},\times))$.
For this expression to make sense (even as a formal power series) we must assume that the numbers $\#(X(\F_{p^n}))$ are all finite.

\begin{proposition}
The formal power series $Z_X(p,T)$ defined above coincides with the Hasse-Weil zeta function of $X\otimes_{\F_1}\F_p=X_{\F_p}$.
\end{proposition}

\prf
This is an immediate consequence of Theorem \ref{adjoint}.
\qed

This type of zeta function thus does not give new insights.
Recall that to get a zeta function over $\Z$, one considers the product
\begin{eqnarray*}
Z_{X\otimes\Z}(s) &=& \prod_{p} Z_{X_{\F_p}}(p^{-s})\\
&=& \prod_p Z_X (p,p^{-s}).
\end{eqnarray*}
As this product takes care of the fact that the prime numbers are the prime places of $\Z$, over $\F_1$ there is only one place, so there should be only one such factor.
Soul\'e \cite{Soule}, inspired by Manin \cite{Manin}, provided  the following idea: Suppose there exists a polynomial $N(x)\in\Z[x]$ such that, for every $p$ one has $\# X(\F_{p^n})=N(p^n)$. 
Then $Z_X(p,p^{-s})$ is a rational function in $p$ and $p^{-s}$. 
One can then ask for the value of that function at $p=1$.
The (vanishing-) order at $p=1$ of $Z_X(p,p^{-s})^{-1}$ is $N(1)$, so the following limit exists
$$
\zeta_X(s)\df \lim_{p\to 1}\frac{Z_X(p,p^{-s})^{-1}}{(p-1)^{N(1)}}.
$$
One computes that if $N(x)=a_0+a_1x+\dots +a_n x^n$, then
$$
\zeta_X(s)\= s^{a_0} (s-1)^{a_1}\cdots (s-n)^{a_n}.
$$
For example $X=\Spec\F_1$ gives
$$
\zeta_{\Spec\F_1}(s)\= s.
$$
For the affine line $\A_1=\Spec C_{\infty,+}$  one gets $N(x)=x$ and thus
$$
\zeta_{\A_1}(s)\= s-1.
$$
Finally, for $\GL_1$ one gets
$$
\zeta_{\GL_1}(s)\= \frac s{s-1}.
$$
In our context the question must be if we can retrieve these zeta functions from the monoidal viewpoint without regress to the finite fields $\F_{p^n}$?
In the examples it indeed turns out that
$$
N(k)\= \# X(\F_{D_k}),
$$
where $D_k$ is the monoid $C_{k-1}\cup \{ 0\}$ and $C_{k-1}$ is the cyclic group with $k-1$ elements.
Since $(\F_{p^n},\times)\cong D_{p^n}$ this comes down to the following question.

\vspace{10pt}

\noindent
{\bf Question.}
Let $X$ be a scheme over $\F_1$.
Assume there is a polynomial $N(x)$ with integer coefficients such that $\# X(D_{p^n})=N(p^n)$ for every prime number $p$ and every non-negative integer $n$.
Is it true that $\# X(D_k)=N(k)$ for every $k\in\N$?
Or, another question: is there a natural characterization of the class of schemes $X$ over $\F_1$ for which there exists a polynomial $N_X$ with integer coefficients such that $\# X(D_k)=N_X(k)$ for every $k\in\N$?

\noindent
{\small Mathematisches Institut, Auf der Morgenstelle 10, 72076 T\"ubingen, Germany, \tt deitmar@uni-tuebingen.de}

\end{document}